\documentclass{amsart}
\usepackage[margin=1.5in]{geometry}

\usepackage{amsfonts, amsmath, amssymb, amsthm}
\usepackage[mathscr]{eucal} 
\usepackage{paralist}
\usepackage{hyperref}

\usepackage{graphicx,color}
\usepackage{booktabs}
\usepackage{bbm}
\usepackage{xspace}

\usepackage{tikz}
\usetikzlibrary{fit,backgrounds,shapes}

\newcommand{\norm}[1]{\left\lVert#1\right\rVert}

\newcommand{\NN}{{\mathbb{N}}}

\newcommand{\ZZ}{{\mathbb{Z}}}
\newcommand{\RR}{{\mathbb{R}}}

\newcommand{\cA}{{\mathscr{A}}}

\newcommand{\cH}{{\mathscr{H}}}

\newcommand{\vones}{\mathbbm{1}}

\newcommand{\SetOf}[2]{\left\{#1\vphantom{#2}\,\right.\left|\,\vphantom{#1}#2\right\}}
\newcommand{\smallSetOf}[2]{\{#1\,|\,#2\}}

\newcommand{\order}[1]{O(#1)}

\DeclareMathOperator{\lin}{lin}

\DeclareMathOperator{\conv}{conv}
\DeclareMathOperator{\pos}{pos}

\DeclareMathOperator{\Sym}{Sym}

\DeclareMathOperator{\id}{id}

\newcommand{\smallpmatrix}[1]{\left(\begin{smallmatrix} #1 \end{smallmatrix}\right)}

\newcommand\Aut{\operatorname{Aut}} 
\newcommand\GL{{\operatorname{GL}}} 
\newcommand\Orth{{\operatorname{O}}} 
\newcommand\Fix{\operatorname{Fix}} 
\newcommand\LP{\operatorname{LP}} 
\newcommand\ILP{\operatorname{ILP}} 
\newcommand\ILPgraph{\operatorname{G}} 
\newcommand\feasible{P} 
\newcommand\arrangement{\cH} 
\newcommand\colvectwo[2]{\begin{pmatrix}{#1}\\ {#2}\end{pmatrix}} 

\newcommand{\Eig}{\operatorname{Eig}} 
\newcommand\coprime[1]{\operatorname{coprime}(#1)} 

\theoremstyle{plain}
\newtheorem{theorem}{Theorem}
\newtheorem{proposition}[theorem]{Proposition}
\newtheorem{corollary}[theorem]{Corollary}
\newtheorem{lemma}[theorem]{Lemma}

\theoremstyle{definition}
\newtheorem{definition}[theorem]{Definition}

\newtheorem{example}[theorem]{Example}

\newtheorem{remark}[theorem]{Remark}

\usepackage[boxed,vlined,oldcommands]{algorithm2e}
\SetKwInOut{Input}{input}
\SetKwInOut{Output}{output}
\SetArgSty{rm}
\SetFuncSty{tt}
\setalcapskip{1ex}
 
\renewcommand\thealgocf{\Alph{algocf}}

\begin{document}
\setdefaultitem{$\triangleright$}{}{}{}

\title[Symmetric Linear and Integer Programs]{Algorithms for Highly Symmetric\\ Linear and Integer Programs}

\author{Richard B\"odi}
\address{School of Engineering, Z\"urcher Hochschule f\"ur Angewandte Wissenschaften, 
Rosenstr.~2, 8400 Winterthur, Switzerland}
\email{richard.boedi@zhaw.ch}
\author{Katrin Herr}
\address{Fachbereich Mathematik, TU Darmstadt, Dolivostr.~15, 64293 Darmstadt, Germany}
\email{\{herr,joswig\}@mathematik.tu-darmstadt.de}
\author{Michael Joswig}
\thanks{Research by Herr is supported by Studienstiftung des deutschen Volkes. Research by
  Joswig is partially supported by DFG Research Unit 565 ``Polyhedral Surfaces'' and DFG
  Priority Program 1489 ``Experimental Methods in Algebra, Geometry, and Number Theory''.}

\date{\today}

\begin{abstract}
  This paper deals with exploiting symmetry for solving linear and integer programming
  problems.  Basic properties of linear representations of finite groups can be used to
  reduce symmetric linear programming to solving linear programs of lower dimension.
  Combining this approach with knowledge of the geometry of feasible integer solutions
  yields an algorithm for solving highly symmetric integer linear programs which only
  takes time which is linear in the number of constraints and quadratic in the dimension.
\end{abstract}

\maketitle

\section{Introduction}
\noindent
It is a known fact that many standard (integer) linear programming formulations of
relevant problems in optimization show a lot of symmetry.  In this situation a standard
branch-and-cut framework repeatedly enumerates symmetric solutions, and sometimes this
renders such methods useless. To address these issues the last decade saw a number of
approaches to devise algorithms specialized to symmetric optimization problems.  We
mention a few: Margot suggests to solve symmetric integer linear programs (ILPs) via a
pruned branch-and-cut approach involving techniques from computational group
theory~\cite{Margot02,Margot03}.  A recent improvement in this direction is ``orbital
branching'' devised by Ostrowski et al.~\cite{OstrowskiEtAl09}.
Friedman~\cite{Friedman07} as well as Kaibel and Pfetsch~\cite{KaibelPfetsch08} treat
symmetric ILPs by shrinking the domain of feasibility by cutting off symmetric
solutions. Gatermann and Parrilo apply results from representation theory and invariant
theory to semidefinite programming~\cite{GatermannParrilo04}, which includes linear
programming as a special case.  Our approach is close in spirit to this paper.  See also
the survey of Margot~\cite{Margot10} for a general overview of symmetric integer linear
programming.

This is how our paper is organized: first we analyze linear programs with an arbitrary
finite group of linear automorphisms.  Most results in this section are known.  A first
key observation, Theorem~\ref{thm:barycenter_ineqs}, is that symmetric linear programming
can be reduced to linear programming over the fixed space of the automorphism group.
Sections~\ref{sec:ILP} and~\ref{sec:layers} translate these results to the context of
integer linear programming.  In the sequel we concentrate on groups acting as signed
permutations on the standard basis of $\RR^n$.  Section~\ref{sec:searching} contains our
main contribution: our \emph{Core Point Algorithm}~\ref{algo:core_point} can solve an
integer linear program in $\RR^n$ whose group of linear automorphisms contains the
alternating group of degree~$n$ (acting as signed permutations) in $\order{mn^2}$ time,
where $m$ is the number of constraints.  This is in sharp contrast with the known
NP-completeness of the general integer linear programming feasibility problem.  While our
algorithm only works for ILPs with an exceptionally high degree of symmetry we believe
that this is a first step towards an entire new class of algorithms dealing with symmetry
in ILPs.  Suitable generalizations are the subject to ongoing research.  In
Section~\ref{sec:finding} we discuss algorithms to determine groups of automorphisms of
integer linear programs.  This leads to algorithmic problems which turn out to be
graph-isomorphism-complete.  The final Section~\ref{sec:experiments} contains experimental
results.  One of the ILP classes that we investigated computationally is motivated by work
of Pokutta and Stauffer on lower bounds for Gomory-Chv\'atal
ranks~\cite{PokuttaStauffer10}.  Section~\ref{sec:hypertruncated} explains the
construction.

We are indebted to Tobias Achterberg, Leo Liberti, Marc Pfetsch, Sebastian Pokutta, and
Achill Sch\"urmann for valuable discussions on the subject of this paper.

\section{Automorphisms of Linear Programs}
\label{sec:LP}
\noindent
The purpose of this section is to introduce the notation and to collect basic facts for
future reference.  The results of this section up to and including
Corollary~\ref{cor:opt_sol_in_fix} can be found in the literature which is why we skip
some of the proofs.

We consider linear programs $\LP(A,b,c)$ of the form
\begin{equation}\label{eq:LP}
  \begin{array}{ll}
    \max & c^tx\\
    \textrm{s.t.} & Ax \le b \, , \ x \in \RR^n
  \end{array}
\end{equation}
where $A\in\RR^{m\times n}$, $b\in\RR^m$, and $c\in\RR^n\setminus\{0\}$. Throughout we
will assume that the set $\feasible(A,b):=\smallSetOf{x\in\RR^n}{Ax\le b}$ of
\emph{feasible points} is not empty, and hence it is a convex polyhedron, which may be
bounded or unbounded.  We will also assume that an optimal solution exists.  This is to
say, our linear program $\LP(A,b,c)$ is bounded even if the feasible region may be
unbounded.  In this case the set of optimal solutions forms a non-empty face of
$\feasible(A,b)$.  Our final assumption $c\ne 0$ is not essential for the algorithms
below, but it allows to simplify the exposition somewhat.

Each row of the matrix $A$ corresponds to one linear inequality.  Suppose that one of
these rows is the zero vector.  Then the corresponding right hand side must be
non-negative, since otherwise the linear program would be infeasible, and this was
explicitly excluded above.  But then this linear inequality is trivially satisfied.
Therefore we will further assume that the matrix $A$ does not contain any zero rows.  In
this case each row defines an affine hyperplane.  This way $\LP(A,b,c)$ gives rise to an
arrangement $\arrangement(A,b)$ of $m$ labeled affine hyperplanes in $\RR^n$.

\begin{definition}\label{def:automorphism}
  An \emph{automorphism} of the linear program $\LP(A,b,c)$ is a linear transformation in
  $\GL_n\RR$ which induces a permutation of $\arrangement(A,b)$, which leaves
  $\feasible(A,b)$ invariant, and which does not change the objective value $c^tx$ for any
  feasible point $x\in\feasible(A,b)$.
\end{definition}

The objective function is linear, and hence it follows that an automorphism of
$\LP(A,b,c)$ does not change the objective value on the linear span $\lin(P(A,b))$ of the
feasible points.  One could also take more general affine transformations into account.
In all what comes below this would require a number of straightforward changes.  We
refrain from doing so for the sake of a clearer exposition.  The following examples show
that the three properties to be satisfied by a linear automorphism are mutually
independent.

\begin{example}
  For $m=n=1$ let $A=1$, $b=0$, and $c=1$.  The feasible region is the non-positive ray
  $\RR_{\le 0}$.  Multiplication with any positive real number $\gamma$ leaves the
  feasible region and the hyperplane arrangement (consisting of the origin) invariant. If
  $\gamma\ne 1$ the objective function is not preserved.
\end{example}

\begin{example}
  For $m=n=2$ let
  \[
  A \ = \ \begin{pmatrix} -1&0\\ 0&-1 \end{pmatrix} \, , \quad b \ = \ 0 \, , \quad c \ = \
  \colvectwo{-1}{0} \, .
  \]
  Then $\feasible(A,b)$ is the non-negative quadrant in $\RR^2$.  Now $\gamma =
  \smallpmatrix{1&0 \\ 0&-1}$ leaves the coordinate hyperplane arrangement
  $\arrangement(A,b)$ invariant, but it changes the feasible region.  For each $x\in\RR^2$
  we have $c^tx=c^t\gamma x$.
\end{example}

\begin{example}\label{exmp:redundant}
  For $m=3$ and $n=2$ let
  \[
  A \ = \ \begin{pmatrix} -1 & 0\\ 0 & -1\\ -1 & -2 \end{pmatrix} \, , \quad b \ = \ 0 \, , \quad
  c \ = -\vones \, .
  \]
  The feasible region is the non-negative quadrant in $\RR^2$; the third inequality is
  redundant.  The linear transformation $\gamma =\smallpmatrix{0 & 1 \\ 1 & 0}$ leaves the
  feasible region invariant, and it satisfies $c^tx=c^t\gamma x$ for all $x\in\RR$.  The
  hyperplane arrangement $\arrangement(A,b)$ is changed.
\end{example}

For more examples see \cite{BoediGrundhoeferHerr10}.  There it is also shown that each
finite (permutation) group occurs as the group of automorphisms of a linear program.

\begin{remark}\label{rem:scale}
  It is always possible to scale the rows of the extended matrix $(A|b)$ such that the
  leftmost non-zero coefficient is $\pm1$.  This allows to remove duplicate inequalities
  from the input by sorting.  The complexity of sorting the rows by pair-wise comparison
  is of order $\order{mn\log m}$.  This can be neglected in the asymptotic analysis of our
  algorithms below since it is always dominated.  This way we can always assume that the
  hyperplanes in $\arrangement(A,b)$, that is, the inequalities, and the rows of the
  extended matrix $(A|b)$ are in a one-to-one correspondence.  In the rational case it is
  more natural to scale the inequalities to integer coefficients which are coprime.  This
  is what we will usually do.  For a more sophisticated algorithm to sort out equivalent
  constraints, see Bixby and Wagner~\cite{BixbyWagner87}.
\end{remark}

Since we view points in $\RR^n$ as column vectors, a matrix $\gamma$ representing a
linear transformation acts by multiplication on the left.  The adjoint action on the row
space, and thus on the set of linear inequalities, is by multiplication of the inverse
transpose $\gamma^{-t}$ on the right.  The set of linear transformations permuting the
arrangement $\arrangement(A,b)$ forms a closed subgroup of $\GL_n\RR$, that is, a linear
Lie group.  Similarly, the set of linear transformations leaving the feasible region
$\feasible(A,b)$ invariant forms a linear Lie group.  It follows that the set
$\Aut(\LP(a,b,c))$ of automorphisms of the linear program $\LP(A,b,c)$ also forms a linear
Lie group.  For basic facts about (linear) Lie groups, see Rossmann~\cite{Rossmann02}.

\begin{remark}
  Clearly, the value and the set of optimal solutions of a linear program only depend on
  the non-redundant constraints.  At the expense of one linear program per constraint one
  can get rid of the redundant ones.  This obviously does not help to reduce the
  complexity of solving the linear program given since the linear program for a redundancy
  check is of the same size.  However, for more costly algorithmic problems, like integer
  programming as is discussed below, this reduction can be useful.  In particular, this
  will be the case when the group of automorphisms becomes larger, see
  Example~\ref{exmp:redundant}.  Notice that the notion of ``invariance''
  from~\cite[Definition~3.1]{GatermannParrilo04}, specialized to linear programming,
  implies that redundant constraints are going to be ignored.
\end{remark}

The following result is a consequence of the combinatorial properties of a convex polytope
$P$: the faces of $P$ are partially ordered by inclusion, and this partially ordered set
forms a lattice.  The automorphisms of this lattice, the \emph{face lattice} of $P$, are
called \emph{combinatorial automorphisms}.  Each linear (or affine or projective)
automorphism of $P$ induces a combinatorial automorphism, but, in general, a polytope may
have many combinatorial automorphisms which are not linearly induced.  See
Ziegler~\cite{Ziegler95} for the details.

\begin{lemma}\label{lem:Aut_finite}
  If the feasible region $\feasible(A,b)$ is bounded and full-dimensional, then the
  automorphism group $\Aut(\LP(A,b,c))$ is finite.  Moreover, the objective function $c$
  satisfies $c^t\gamma x=c^t x$ for all $x\in\RR^n$ and $\gamma\in\Aut(\LP(A,b,c))$.
\end{lemma}

\begin{proof}
  Let $v$ be a vertex of the polytope $P=\feasible(A,b)$.  Since $\dim P=n$ there are
  vertices $w_1,w_2,\dots,w_n$ each of which shares an edge with $v$ and such that the
  difference vectors $w_1-v,w_2-v,\dots,w_n-v$ form a basis of $\RR^n$.  This implies that
  each combinatorial automorphism of $P$ is induced by at most one linear automorphism.
  Hence the group $\Aut(\LP(A,b,c))$ is contained in the group of combinatorial
  automorphisms of $P$, which is finite.  While Definition~\ref{def:automorphism} asks
  that each \emph{feasible} point is mapped to a (feasible) point with the same objective
  value, the additional claim deals with all points, feasible or not.  However, this
  follows from $\lin(P(A,b))=\RR^n$ as $c$ is linear.
\end{proof}

If the polyhedron $\feasible(A,b)$ is not full-dimensional, then the automorphism group is
a direct product of the group of automorphisms fixing the linear span of $\feasible(A,b)$
with a full general linear group of the orthogonal complement.  In the sequel we will
therefore restrict our attention to the full-dimensional case.

\begin{definition}
  Given a subset $Y\subseteq \RR^n$ and a group $\Gamma\le\GL_\RR$ acting on $Y$, the
  \emph{set of fixed points} of $Y$ with respect to an element $\gamma\in\Gamma$ is
  defined by
  \[
  \Fix_\gamma(Y) \ := \ \SetOf{y\in Y}{\gamma y=y} \, .
  \]
  Therefore, the \emph{set of fixed points} of $Y$ with respect to $\Gamma$ is given by
  \[
  \Fix_\Gamma(Y) \ := \ \SetOf{y\in Y}{\gamma y=y \text{ for all } \gamma\in\Gamma} \ = \
  \bigcap_{\gamma\in\Gamma}\Fix_\gamma(Y) \, .
  \]
\end{definition}

The set of fixed points $\Fix_\gamma(\RR^n)$ is the (possibly zero-dimensional) eigenspace
$\Eig(\gamma,1)$ of the linear transformation $\gamma$ with respect to the eigenvalue~$1$.
This implies that $\Fix_\Gamma(\RR^n)$ is a linear subspace for any group $\Gamma$ of
linear transformations.  More generally, $\Fix_\Gamma(Y)$ is the intersection of this
subspace with the set $Y$.

\begin{remark}\label{rem:compute_fixed}
  If the linear group $\Gamma\le\GL_n\RR$ is generated by the set $G\subset\Gamma$, then
  \[
  \Fix_\Gamma(\RR^n) \ = \ \bigcap_{\gamma\in G}\Fix_\gamma(\RR^n) \ = \
  \bigcap_{\gamma\in G}\Eig(\gamma,1) \, .
  \]
  In particular, if $G$ is finite, that is, if the group $\Gamma$ is finitely generated,
  this leads to an algorithm to compute (a primal or dual basis of) the fixed space by
  solving one linear system of equations per transformation in the
  generating set~$G$.
\end{remark}  

\begin{remark}\label{rem:c_one}
  Let $\Gamma\le\Aut(\LP(A,b,c))$ be a group of automorphisms of the linear program
  $\LP(A,b,c)$ such that $P(A,b)$ is bounded and full-dimensional.  Then, by
  Lemma~\ref{lem:Aut_finite}, the set of fixed points $\Fix_\Gamma(\RR^n)$ contains the
  one-dimensional linear subspace spanned by the objective vector~$c$.
\end{remark}

For any finite set $S\subset \RR^n$ let
\[
\beta(S) \ := \ \frac{1}{|S|}\sum_{v\in S}v
\]
be its \emph{barycenter}.  The two subsequent results are basic observations from
representation theory, not restricted to applications in (linear) optimization.  For a brief
proof, for instance, see~\cite[Lemma~3.5]{Webb}.

\begin{lemma}\label{lem:barycenter}
  The map
  \[
  \RR^n \to \Fix_\Gamma(\RR^n) \, , \ v \mapsto \beta(\Gamma v)
  \]
  which sends a point to the barycenter of its $\Gamma$-orbit is a linear projection onto
  the fixed space.
\end{lemma}


Let $S\subseteq \RR^n$ be a finite set which is \emph{spanning}, that is, we require
$\lin(S)=\RR^n$.  Further let $\Gamma$ be a finite subgroup of $\GL_n\RR$ acting on $S$.
Phrased differently, we are considering a linear representation of an abstract group
$\Gamma$ on the vector space $\RR^n$ which induces a permutation representation on the set
$S$.  In this case $\Gamma$ splits $S$ into disjoint orbits $O_1,O_2,\ldots,O_k$.  In our
applications below, $S$ will usually be the set of vertices of some polytope which
linearly spans $\RR^n$.

\begin{lemma}\label{lem:fix_spanned_by_betaO}
  For the fixed space of $\Gamma$ we have
  \[
  \Fix_\Gamma(\RR^n) \ = \ \lin\{\beta(O_1),\beta(O_2),\dots,\beta(O_k)\} \, .
  \]
  In particular, $\dim\Fix_\Gamma(\RR^n) \le k$.
\end{lemma}

\begin{proof}
  Since $S=O_1\cup O_2\cup\dots\cup O_k$ is spanning and since the union of the orbits
  gives $S$ it follows that
  \begin{equation}\label{eq:orbit_decomp}
    \RR^n \ = \ \lin(O_1)+ \lin(O_2)+ \dots+ \lin(O_k) \, .
  \end{equation}
  For $i\in\{1,2,\dots,k\}$ the linear subspace $\lin(O_i)$ is $\Gamma$-invariant.  If
  we apply the surjective linear map $v\mapsto\beta(\Gamma v)$ from
  Lemma~\ref{lem:barycenter} to the set $S$, we obtain a generating set for
  $\Fix_\Gamma(\RR^n)$.  Applying the same map to a single orbit $O_i$ similarly yields a
  generating set for $\Fix_\Gamma(\lin(O_i))$.  Now the claim follows from the
  equation $\Gamma O_i=O_i$.
\end{proof}

Notice that the sum decomposition \eqref{eq:orbit_decomp} is not necessarily direct.  We
now apply the results obtained so far to a finite group of automorphisms of a linear
program.

\begin{proposition}\label{prop:symmetric_solution}
  Let $\Gamma\le\Aut(\LP(A,b,c))$ be finite.  If $x\in\RR^n$ is an arbitrary point, the
  barycenter of its $\Gamma$-orbit satisfies $c^t \beta(\Gamma x)=c^t x$.  If, moreover,
  $x\in \feasible(A,b)$ is feasible, then $\beta(\Gamma x)$ is feasible, too.
\end{proposition}

Geometrically this means that the points of one orbit are in the same affine hyperplane
orthogonal to $c$.

\begin{proof}
  As the objective function is constant on the orbit $\Gamma x$ it follows that $c^t
  \beta(\Gamma x)=c^t x$.  If $x$ is a feasible point, then $\gamma x$ is also feasible for
  all $\gamma\in\Gamma$.  So the barycenter $\beta(\Gamma x)$ is a convex combination of
  feasible points.  The claim follows as the feasible region is convex.
\end{proof}

Since we assumed that $\LP(A,b,c)$ has an optimal solution, the following is an immediate
consequence of the preceding result.

\begin{corollary}\label{cor:opt_sol_in_fix}
  There exists an optimal solution of $\LP(A,b,c)$ which is a fixed point with respect to
  the entire automorphism group $\Aut(\LP(A,b,c))$.
\end{corollary}

Up to minor technical details Theorem~3.3 of~\cite{GatermannParrilo04} generalizes
Corollary~\ref{cor:opt_sol_in_fix} to semi-definite programming.

Let $\LP(A,b,c)$ be a linear program with $\feasible(A,b)$ bounded and full-dimensional,
and let $\Gamma=\langle \gamma_1,\gamma_2,\dots,\gamma_t\rangle$ be a finite subgroup of
$\Aut(\LP(A,b,c))$.  Following Remark~\ref{rem:compute_fixed} we can compute a matrix $E$
such that the kernel $\smallSetOf{x}{Ex=0}$ is the fixed space $\Fix_\Gamma(\RR^n)$: for
each $\gamma_i$ we determine a dual basis for the eigenspace
$\smallSetOf{x}{(\gamma_i-\id)x=0}$ by solving a square system of linear equations.  The
total number of operations to do so is of order $\order{tn^3}$.  Throughout this paper we
measure algorithmic complexity in the RAM model; that is, we ignore the encoding lengths
of real numbers, and all arithmetic operations are assumed to take constant time.  The
group $\Gamma$ acts on the rows of the extended matrix $(A|b)$, and we define a new
extended matrix $(A'|b')$ by summing the rows of the same $\Gamma$-orbit.  We have the
following general result.

\enlargethispage{2em}

\begin{theorem}\label{thm:barycenter_ineqs}
  The polyhedron
  \[
  P' \ = \ \SetOf{x\in\RR^n}{A'x\le b',\, Ex=0}
  \]
  is the set $\Fix_\Gamma(P(A,b))$ of feasible points which is fixed under the action of
  $\Gamma$.  In particular, $P'=\smallSetOf{\beta(\Gamma x)}{x\in P(A,b)}$.  Each optimal
  solution of the linear program
  \begin{equation}\label{eq:reduced_LP}
    \begin{array}{ll}
      \max & c^tx\\
      \textrm{s.t.} & \begin{pmatrix} A'\\ E\\ -E\end{pmatrix} x \le \begin{pmatrix}b'\\ 0\\
        0\end{pmatrix} \, , \ x \in \RR^n
    \end{array}
  \end{equation}
  is an optimal solution of $\LP(A,b,c)$, and the objective values are the same.
\end{theorem}

\begin{proof}
  We constructed the matrix $E$ to guarantee that each fixed point in $P=P(A,b)$ satisfies
  the equation $Ex=0$.  Further, each inequality of the system $A'x\le b'$ is a positive
  linear combination of valid inequalities.  It follows that $\Fix_\Gamma(P)$ is contained
  in $P'$.

  To prove the reverse inclusion consider a point $x$ which is fixed by each
  transformation in $\Gamma$ but which is not contained in $P$.  Then for some index $i$
  we have the strict inequality $a_{i,\cdot}x>b_i$.  Without loss of generality we can
  assume that the first $k$ rows $a_{1,\cdot},a_{2,\cdot},\dots,a_{k,\cdot}$ of $A$ form
  the $\Gamma$-orbit of the row $a_{i,\cdot}$.  It follows that $b_1=b_2=\dots=b_k=b_i$.
  Moreover, since $x$ is a fixed point we have
  \[
  a_{1,\cdot}x \ = \ a_{2,\cdot}x \ = \ \cdots \ = \ a_{1,\cdot}x \ = a_{i,\cdot}x \ > \
  b_i \, .
  \]
  This implies that $(\sum_{j=1}^k a_{j,\cdot})x>kb_i$, and hence $x$ is not contained in
  $P'$.  We conclude that $P'$ is the set of points in $P$ fixed by each transformation of
  $\Gamma$. Now Lemma~\ref{lem:barycenter} says that $P'$ is the image of $P$ under the
  map $x\mapsto \beta(\Gamma x)$.  The claim about the linear
  program~\eqref{eq:reduced_LP} follows from Corollary~\ref{cor:opt_sol_in_fix}.
\end{proof}

\begin{remark}
  It has been observed by Scharlau and Sch\"urmann\footnote{private communication} that
  the vertices of the polyhedron $P'$ are barycenters of orbits of vertices of $P$.  This
  is a consequence of the fact that $P'$ is the image of $P$ under the linear map
  $x\mapsto\beta(\Gamma x)$.
\end{remark}

Corollary~\ref{cor:opt_sol_in_fix} and Theorem~\ref{thm:barycenter_ineqs} yield a direct
algorithm for solving a symmetric linear program: instead of solving $\LP(A,b,c)$ one can
solve the linear program~\eqref{eq:reduced_LP}.  The benefit is the following: The larger
the group $\Gamma\leq\Aut(\LP(A,b,c))$ the smaller the dimension of the fixed space and
the number of constraints.

\begin{remark}
  Formally, the feasible points of the derived linear program live in the same space
  $\RR^n$ as the original linear program.  However, an algorithm based on the Simplex
  Method directly benefits if the solutions are contained in a proper subspace: the rows
  of the matrix $E$ describing the fixed space never have to be exchanged in a Simplex
  tableau.  Alternatively, one can project $\Fix_\Gamma(\RR^n)$ onto a full-dimensional
  coordinate subspace, solve the projected linear program and lift back.
\end{remark}

In the special case where the linear program admits a group of automorphisms acting on the
standard basis of $\RR^n$ (that is, the groups acts by permuting the columns) it is
standard optimization practice to identify variables in the same orbit, and to solve the
reduced linear program.  Theorem~\ref{thm:barycenter_ineqs} generalizes this approach to
arbitrary groups of automorphisms.

\section{Symmetries of Integer Linear Programs}
\label{sec:ILP}
\noindent
We now turn to our main focus.  Associated with $\LP(A,b,c)$ is the integer linear program
\begin{equation}\label{eq:ILP}
  \begin{array}{ll}
    \max & c^tx\\
    \textrm{s.t.} & Ax \le b \, , \ x \in \ZZ^n  \, ,
  \end{array}
\end{equation}
which we denote as $\ILP(A,b,c)$.  Throughout we make the same assumptions as above: the
linear program $\LP(A,b,c)$ is feasible, the matrix $A$ does not have any zero rows, and
the inequalities bijectively correspond to the hyperplane arrangement $\arrangement(A,b)$;
see Remark~\ref{rem:scale}.

\begin{definition}
  A \emph{symmetry} of the integer linear program $\ILP(A,b,c)$ is an automorphism of
  $\LP(A,b,c)$ which acts on the signed standard basis $\{\pm e_1,\pm e_2,\dots,\pm e_n\}$
  of $\RR^n$ as a signed permutation.
\end{definition}

The symmetries of the integer linear program \eqref{eq:ILP} form a group
$\Sym(\ILP(A,b,c))$ which is a subgroup of the group $\Orth_n\ZZ$, the group of all
$0/1/{-}1$-matrices with exactly one non-zero entry per row and column.  We have
$\Orth_n\ZZ=\Orth_n\RR\cap\GL_n\ZZ$, and $\Orth_n\ZZ$ is isomorphic to the Coxeter group
of type $B_n$, the group of automorphisms of the regular $n$-dimensional cube and its
polar, the regular $n$-dimensional cross polytope.  As a consequence, the group of
symmetries of an integer linear program is finite, even if $\Aut(\LP(A,b,c))$ is infinite.

The motivation for our definition is Lie-theoretical: let $\Gamma$ be any finite subgroup
of $\GL_n\ZZ$.  Then $\Gamma$ is a compact subgroup of $\GL_n\RR$, hence it is contained
in (a conjugate copy of) the maximal compact subgroup $\Orth_n\RR$.  It follows that, up
to conjugation in $\GL_n\RR$, the group $\Gamma$ is a subgroup of $\Orth_n\ZZ$.

As an abstract group $\Orth_n\ZZ$ is isomorphic to the wreath product
\[
\ZZ_2\wr\Sym(n) \ = \ (\ZZ_2)^n \rtimes \Sym(n) \, ,
\]
where $\ZZ_2$ is the cyclic group of order two and $\Sym(n)$ is the symmetric group of
degree $n$; the group $\Sym(n)$ acts on the direct product $(\ZZ_2)^n$ by permuting the
factors.  Each element of $\Orth_n\ZZ$ can be written as a product of a sign vector and a
permutation.  Since a permutation is a product of disjoint cycles, each signed permutation
is a product of signed cycles which are disjoint.  In terms of notation we write the signs
between the indices within a cycle.  This is to say, $(1{+}2{-}4{+}3{-})$ denotes the
signed permutation matrix
\[
\begin{pmatrix}
  0 & 0 & -1& 0\\
  1 & 0 & 0 & 0\\
  0 & 0 & 0 & 1\\
  0 & -1 & 0 & 0
\end{pmatrix}
\]
which is to be multiplied to column vectors from the left.

\begin{remark}
  In the optimization literature the authors often restrict their attention to symmetries
  permuting the standard basis vectors; for instance, see Margot~\cite{Margot10} and the
  references listed there.  However, our more general analysis below shows that taking
  signed permutations into account does not cause any extra effort.  Moreover, if the
  polyhedron $P(A,b)$ is full-dimensional and bounded the group of automorphisms of the
  linear relaxation is already finite by Lemma~\ref{lem:Aut_finite}.  Then
  $\Aut(\LP(A,b,c))\cap\GL_n\ZZ$ is already contained in $\Orth_n\ZZ$ by the
  Lie-theoretical argument given above.  Hence, at least in this case, considering
  groups of signed permutations is a natural choice.
\end{remark}

Before we will inspect groups of symmetries of integer linear programs we need to collect
a few basic results on the action of the group $\Orth_n\ZZ$ on the entire space~$\RR^n$.
Throughout let $\Gamma$ be a subgroup of $\Orth_n\ZZ$.  Then $\Gamma$ acts on the standard
basis
\[
S \ = \ \{\pm e_1,\pm e_2,\dots,\pm e_n\} \, .
\]
In the sequel we will always consider this particular action of $\Gamma$.  There are two
kinds of orbits to distinguish: the \emph{bipolar} orbits contain at least one pair $\pm
e_i$, while the \emph{unipolar} orbits do not.  Since $\Gamma$ is a linear group, a signed
permutation $\sigma\in\Gamma$ with $\sigma e_i=\epsilon e_j$ and $\epsilon\in \{\pm 1\}$
maps $-e_i$ to $-\epsilon e_j$.  Hence, a bipolar orbit only consists of pairs, that is,
$-O=O$.  On the other hand, for each unipolar orbit $O$ the set
$-O=\smallSetOf{-e_i}{e_i\in O}$ forms another orbit, and $\Gamma$ acts equivalently on
$O$ and~$-O$.

\begin{proposition}\label{prop:fix} 
  For the fixed space of $\Gamma$ we have
  \[
  \Fix_\Gamma(\RR^n) \ = \ \lin\SetOf{\beta(O)}{\text{$O$ orbit of $\Gamma$}} \ = \
  \lin\SetOf{\beta(O)}{\text{$O$ unipolar orbit of $\Gamma$}} \, .
  \]
\end{proposition}

\begin{proof}
  The first equality is a consequence of Lemma~\ref{lem:fix_spanned_by_betaO}.  The second
  equality holds as $\beta(O)=0$ for any bipolar orbit~$O$.
\end{proof}

\begin{remark}
  The points in $S$ are the vertices of the regular $n$-dimensional cross polytope.  If
  $O\subset S$ is a unipolar $\Gamma$-orbit, then $\beta(O)$ is the barycenter of the
  non-trivial face of the cross polytope which is spanned by the vertices in $O$.  In view
  of cone polarity the action of $\Gamma$ on $S$ is dual to the induced action on the
  vertices of the regular cube $[-1,1]$.  That is, the two corresponding representations
  of $\Gamma$, on $\RR^n$ and its dual space, form a contra-gradient pair.
\end{remark}

We call the action of $\Gamma$ on the set of signed standard basis vectors
\emph{semi-transitive} if there are precisely two opposite orbits of length $n$.
Moreover, we call the action \emph{sub-transitive} if there is no proper coordinate
subspace which is invariant under~$\Gamma$.  Clearly, a semi-transitive action is
necessarily sub-transitive.  The converse does not hold, but we have the following
characterization.

\begin{proposition}\label{prop:semi-trs}
  Suppose that $\Gamma$ acts sub-transitively.  Then either $\Gamma$ acts
  semi-transitively with orbits $O$ and $-O$ such that the fixed space
  \[
  \Fix_\Gamma(\RR^n) \ = \ \lin(\beta_O) \ = \ \lin(\beta_{-O})
  \]
  is one-dimensional, or $\Fix_\Gamma(\RR^n)=0$.
\end{proposition}

\begin{proof}
  If $\Gamma$ has a bipolar orbit $O$, then $O$ equals the entire set $S$ of signed
  standard basis vectors because $\Gamma$ acts sub-transitively.  In this case the fixed
  space reduces to the origin.  If, however, each orbit is unipolar, we have exactly one
  pair $(O,-O)$ of opposite orbits, again due to sub-transitivity.  Now the claim follows
  from Proposition~\ref{prop:fix}.
\end{proof}

\begin{corollary}\label{cor:conjugate}
  If $\Gamma$ acts semi-transitively, then $\Gamma$ is conjugate to a subgroup of $\Sym(n)$
  in $\Orth_n\ZZ$.
\end{corollary}

\begin{proof}
  Let $O$ and $-O$ be the two orbits of $\Gamma$, both of which have length~$n$.  Pick a
  transformation $\epsilon\in\Orth_n\ZZ$ which maps $O$ to the standard basis
  $\{e_1,e_2,\dots,e_n\}$.  Now for each $\gamma\in\Gamma$ the conjugate transformation
  $\epsilon\gamma\epsilon^{-1}$ leaves the sets $\{e_1,e_2,\dots,e_n\}$ and
  $\{-e_1,-e_2,\dots,-e_n\}$ invariant.  We conclude that $\epsilon\Gamma\epsilon^{-1}$ is
  a subgroup of $\Sym(n)$.
\end{proof}

We now interprete the results above for integer linear programming.  Consider an integer
linear program $\ILP(A,b,c)$ such that the set $P(A,b)$ of feasible points of the linear
relaxation is full-dimensional.  Let $\Gamma\le\Sym(\ILP(A,b,c))$ be a group of
automorphisms.  We have $\Gamma\le\Orth_n\ZZ$.  The action of $\Gamma$ on the set $\{\pm
e_1,\pm e_2,\dots,\pm e_n\}$ can be decomposed into orbits.  In this way the most relevant
case occurs when $\Gamma$ acts sub-transitively.  From Lemma~\ref{lem:Aut_finite} we know
that $c$ is contained in the fixed space $\Fix_\Gamma(\RR^n)$, and then
Proposition~\ref{prop:semi-trs} says $c\ne 0$ enforces the action of $\Gamma$ to be
semi-transitive.  Finally, by Corollary~\ref{cor:conjugate} we can conjugate $\Gamma$ into
a subgroup of $\Sym(n)$ acting on the standard basis $\{e_1,e_2,\dots,e_n\}$.  This is the
situation that we will be dealing with in our algorithms below.

\section{Layers of Integer Points}
\label{sec:layers}
\noindent
Our goal is to describe an algorithm for the efficient solution of a highly symmetric
integer linear program.  Again we consider $\ILP(A,b,c)$ with a group $\Gamma$ of
automorphisms as above.

Let us assume that the objective function $c\ne 0$ is \emph{projectively rational}.  This
means that we require $c$ to be a constant real multiple of a rational vector.  For such a
vector $c$ let $\coprime{c}$ be the unique integral vector with coprime coefficients such
that $c=\rho\coprime{c}$ for some positive real $\rho$.  If $c$ is a multiple of a
standard basis vector, the single non-zero coefficient of $\coprime{c}$ is defined to be
$\pm 1$.  For an integer $k$ the \emph{$k$-th $c$-layer} is the affine hyperplane
\[
H_{c,k} \ = \ \ker\ (x\mapsto c^tx)+\frac{k}{\norm{\coprime{c}}^2}\coprime{c} \, .
\]
We have $H_{c,k}=H_{\rho c,k}$ for all $\rho>0$, and $H_{-c,k}=-H_{c,k}=H_{c,-k}$.  All
points in $H_{c,k}$ attain the same value $k$ with respect to the rescaled objective
function $\coprime{c}$.  We call $k$ the \emph{number} of the $c$-layer $H_{c,k}$.  The
intersection of $H_{c,k}$ with the line $\RR c$ is called the \emph{center}.

\begin{lemma}\label{lem:layer}
  If $c\ne 0$ is projectively rational, the integral point $x\in\ZZ^n$ is contained in the
  $c$-layer with number $\coprime{c}^tx$.
\end{lemma}

\begin{proof}
  The number $k=\coprime{c}^tx$ is an integer.  We abbreviate $d=\coprime{c}$ and compute 
  \[
  c^t\bigl(\frac{k}{\norm{d}^2}d\bigr) \ = \ c^t\bigl(\frac{d^txd}{\norm{d}^2}\bigr) \ = \ \frac{d^td}{\norm{d}^2}c^tx \ = \ c^tx \, .
  \]
  Hence $x-(k/\norm{d}^2)d$ is contained in the kernel of the linear form $c^t$, that is,
  the point $x$ lies in the affine hyperplane $H_{c,k}$.
\end{proof}

For the following result it is crucial that the coefficients of $\coprime{c}$ are coprime.

\begin{proposition}\label{prop:partition}
  If $c\ne 0$ is projectively rational, the $c$-layers $H_{c,k}$ for $k\in\ZZ$ partition
  the set $\ZZ^n$ of all integral points.
\end{proposition}

\begin{proof}
  From Lemma~\ref{lem:layer} is clear that each integral point is contained in some
  $c$-layer.  By construction it is also obvious that the $c$-layers are pairwise
  disjoint.  It remains to show that $H_{c,k}\cap\ZZ^n$ is non-empty for all $k\in\ZZ$.

  Let $d=\coprime{c}$.  Since the coefficients $d_1,d_2,\dots,d_n$ are coprime there are
  integral coefficients $x_1,x_2,\dots,x_n$ such that
  \[
  x_1 d_1 + x_2 d_2 + \dots + x_n d_n \ = \ \gcd(d_1,d_2,\dots,d_n) \ = \ 1 \, .
  \]
  However, the left side of this equation equals $c^t x$, whence the point $x$ is
  contained in the first $c$-layer $H_{c,1}$.  Now $c^t(kx)=kc^t(x)=k$ implies that the
  $k$-th layer contains the integral point~$kx$ for arbitrary $k\in\ZZ$.
\end{proof}

Another way of putting the statement above is that $\coprime{c}$ is the unique generator
of the unique minimal Hilbert basis of the one-dimensional pointed cone $\RR_{\ge 0}c$.

\begin{remark}
  An important consequence of Proposition~\ref{prop:partition} is that for any given bounds
  $\ell,u\in\RR$ there are only finitely many $c$-layers with feasible integral points whose
  objective values lie between $\ell$ and $u$.  This does not hold if the objective function
  is not projectively rational.
\end{remark}

\begin{algorithm}[tb]
  \dontprintsemicolon
  \linesnumbered

  \KwIn{$(A,b)$ such that $\Sym(\ILP(A,b,\vones))$ acts transitively on
    standard basis}
  \KwOut{optimal solution of $\ILP(A,b,\vones)$ or ``infeasible''}

  let $z=\zeta\vones$ be a symmetric optimal solution of the LP relaxation $\LP(A,b,\vones)$ \;
  $k\leftarrow \lfloor n\zeta\rfloor$ \;
  \Repeat{feasible $x$ found or $k<n\lfloor\zeta\rfloor$}{
    \eIf{it exists}{
      let $x$ be integral point in $P(A,b)\cap H_{\vones,k}$ \; \nllabel{step:reduction_to_layers:feasibility}
    }{
      $k\leftarrow k-1$ \;
    }
  }
  \Return{$x$ or ``infeasible''}
  \caption{Reduction to $\vones$-layers}
  \label{algo:reduction_to_layers}
\end{algorithm}

\begin{theorem}\label{thm:reduction_to_layers}
  For given $A$ and $b$ such that $\Sym(\ILP(A,b,\vones))$ acts transitively on the
  standard basis the Algorithm~\ref{algo:reduction_to_layers} solves the integer linear
  program $\ILP(A,b,\vones)$.
\end{theorem}

\begin{proof}
  Recall that throughout we assumed that the set of feasible points of the linear
  relaxation is bounded.  Hence it cannot occur that the integer linear program is
  unbounded.

  Let $\Gamma\le\Sym(\ILP(A,b,\vones))$ be a transitive group of automorphisms.  The fixed
  space is spanned by $\vones$.  If $z$ is an optimal solution of the relaxation
  $\LP(A,b\vones)$, then, by Proposition~\ref{prop:symmetric_solution}, the barycenter
  $\beta(\Gamma z)=\zeta\vones$ for $\zeta=1/n(z_1+z_2+\dots+z_n)$ is also an optimal
  solution.  Now $\lfloor\zeta\rfloor\vones$ is an integral point in the fixed space with
  an objective value not exceeding the optimal value of the linear programming relaxation.
  Each $\vones$-layer with a feasible integral point meets the one-dimensional polyhedron
  $P'=\smallSetOf{\beta(\Gamma x)}{x\in P(A,b)}$.  We infer that no integral optimal
  solution of $\ILP(A,b,\vones)$ can have an objective value strictly less than
  $n\lfloor\zeta\rfloor$.

  Due to Proposition~\ref{prop:partition} the $\vones$-layers partition $\ZZ^n$, and so
  the feasible points of $\ILP(A,b,c)$ are contained in the set
  \[
  \bigcup_{k=n\lfloor\zeta\rfloor}^{\lfloor n\zeta\rfloor} H_{\vones,k} \, .
  \]
\end{proof}

The benefit of Algorithm~\ref{algo:reduction_to_layers} is that it reduces a (symmetric)
$n$-dimensional integer linear programming problem to $n$ integer feasibility problems in
one dimension below.  Since the latter is still an NP-complete problem not much is gained,
in general.  The situation changes, however, if we assume higher degrees of transitivity
for the action of the group of automorphisms.

\begin{remark}
  Searching a family of parallel affine hyperplanes for integer points as in
  Algorithm~\ref{algo:reduction_to_layers} also plays a key role in Lenstra's algorithm
  for integer linear programming which requires polynomial time in fixed
  dimension~\cite{Lenstra83}.
\end{remark}

\section{Searching Integer Layers Efficiently}
\label{sec:searching}
\noindent
The question remaining is how to test ILP-feasibility of a $c$-layer in an efficient
way. Our key observation is that some optimal integral solution is close to the fixed
space if the group of symmetries acts sufficiently transitive.

\begin{definition}
  Given a $c$-layer with center $z$, an integral point in the $c$-layer is a \emph{core
    point} if it minimizes the distance to $z$.
\end{definition}

\begin{example}\label{exmp:core}
  For the objective function $c=\vones$ and an integer $k=q n+r$ with $q\in\ZZ$ and
  $r\in\{0,1,\dots,n-1\}$, the set of core points in the $k$-th layer consists of all
  integer points with $r$ coefficients equal to $q+1$ and $n-r$ coefficients equal to $q$.
  In particular, the number of core points in this case equals $\tbinom{n}{r}$.  These
  core points are the vertices of an $(r,n)$-hypersimplex, translated by the vector
  $q\vones$.
\end{example}

For the algorithms below the geometric structure of the set of core points is very
relevant.  We therefore make a short digression: the \emph{$(r,n)$-hypersimplex}
$\Delta(r,n)$ is the $0/1$-polytope with vertices
\[
e_S \ = \ \sum_{i\in S} e_i \, ,
\]
where $S$ ranges over all $r$-element subsets of $[n]$.  The hypersimplices are highly
regular structures, and this yields the following.

\begin{proposition}\label{prop:mu-trs}
  Let $\Gamma\le\GL_n\RR$ be a linear group which acts $\mu$-transitively on the standard
  basis.  Then $\Gamma$ acts transitively on the set of vertices of the
  $(r,n)$-hypersimplex for any
  \[
  r \in \{0,1,\dots,\mu\}\cup\{n-\mu,n-\mu+1,\dots,n\} \, .
  \]
\end{proposition}

\begin{proof}
  By assumption $\Gamma$ acts transitively on the $r$-element subsets of $[n]$ for
  $r\le\mu$.  Since $\Gamma$ is a linear group it thus acts transitively on the set of
  vertices of $\Delta(r,n)$.  The corresponding claim for the remaining hypersimplices
  follows since $\Delta(r,n)$ is affinely isomorphic to $\Delta(n-r,n)$ via the map
  $x\mapsto \vones-x$.
\end{proof}

Below we will apply the previous results in the special case where $\mu\ge\lfloor
n/2\rfloor$.  Then the groups acts transitively on the sets of vertices of \emph{all}
hypersimplices.

\begin{lemma}\label{lem:closer}
  Let $x\in\feasible(A,b)$ be an LP-feasible point in the $k$-th $c$-layer, and let
  $\gamma\in\Sym(\ILP(A,b,c))$ with $\gamma x\ne x$.  Then any point in the interior of
  the line segment $[x,\gamma x]$ is LP-feasible and closer to the center of the $k$-th
  $c$-layer than $x$.
\end{lemma}

\begin{proof}
  Since $\gamma$ is an orthogonal linear map it preserves distances.  The center $z$ of
  the $k$-th $c$-layer is fixed under $\gamma$, and this implies that $(x,z,\gamma x)$ is
  an isosceles triangle.  We infer that $\norm{p-z}<\norm{x-z}$ for all points $p$ in the
  interior of $[x,\gamma x]$.  Since $\gamma$ is an automorphism of the linear relaxation
  $\LP(A,b,c)$ the point $\gamma x$ is feasible, too.  The feasible region is convex, and
  hence $p$ is feasible.
\end{proof}

\begin{theorem}\label{thm:core}
  Suppose that $\Gamma\le\Sym(\ILP(A,b,\vones))$ acts $(\lfloor
  n/2\rfloor+1)$-transitively on the standard basis of $\RR^n$, and $n\ge 2$.  Then either
  each core point in the $k$-th $\vones$-layer is feasible or $H_{\vones,k}$ does not
  contain any feasible point.
\end{theorem}

\begin{proof}
  Let $x$ be a feasible integer point in the $k$-th $\vones$-layer which is not a core
  point.  We will show that there is another feasible integer point which is closer to the
  center, and this will prove the claim.

  Due to the invariance of $\ZZ^n$ under translation by integer vectors we may assume that
  $k\in\{0,1,\dots,n-1\}$.  Since $x$ is not a core point, in particular, it is not the
  center of the $k$-th layer.  Hence $x$ is not contained in the fixed space $\RR\vones$,
  which means that not all coordinates of $x$ are the same.  We split the set $[n]$ of
  coordinate directions into two subsets by considering
  \[
  \SetOf{i}{\text{$x_i$ is even}} \quad \text{and} \quad \SetOf{i}{\text{$x_i$ is odd}} \, .
  \]
  Then one of the sets --- denoted by $I$ --- contains at least $\lfloor (n+1)/2\rfloor$
  elements, while the other set $J$ has at most $\lfloor n/2\rfloor$ elements.  We will
  employ the $\lfloor n/2\rfloor$-transitivity of the automorphism group to control $J$,
  and the additional degree of freedom to produce two distinct feasible integer points.
  We distinguish two cases.

  \begin{enumerate}
  \item Suppose that $x$ has two different coordinates, say $x_u$ and $x_v$, which are in
    the same congruence class modulo two.  That is, the set $\{u,v\}$ is contained in
    either $I$ or~$J$.  Observe that this condition is satisfied whenever $x$ has at least
    three pairwise distinct coordinates.  Due to the $(\lfloor n/2\rfloor+1)$-transitivity
    of $\Gamma$ there is an automorphism $\gamma\in\Gamma$ which leaves $J$ invariant and
    which maps $u$ to $v$.  Since $J$ is invariant, its complement $I=[n]\setminus J$ is
    invariant, too.  Notice that we do not require the set $J$ to be non-empty (if
    $\{u,v\}\subseteq I$).

    Letting $x'=\gamma x$ we observe that $x_i$ and $x'_i$ are congruent modulo two for all
    $i\in[n]$.  Since $x_u\ne x_v=x_{\gamma(u)}=x'_u$ we have $x\ne x'$, and hence
    \[
    y \ = \ \frac{1}{2}(x+x') \ = \ \frac{1}{2}(x+\gamma x)
    \]
    is an integer point in the interval $[x,\gamma x]$.
  \item Otherwise the point $x$ has exactly two different coordinates $x_u$ and $x_v$,
    one of them being even, the other one odd.  Without loss of generality, $x_i=x_u$ for
    all $i\in I$ and $x_j=x_v$ for all $j\in J$.  Due to the transitivity of $\Gamma$
    there is an automorphism $\gamma\in\Gamma$ with $\gamma e_u=e_v$.  Then $x$ and
    $\gamma x$ are distinct points.  Consider an interior point
    \begin{equation}
      y \ = \ \lambda x + (1-\lambda)\gamma x \quad \text{for } 0<\lambda<1 \label{eq:y}
    \end{equation}
    in the line segment $[x,\gamma x]$.  We want to find a parameter $\lambda$ such that
    $y$ is integral.  As $x$ has only two distinct coordinates the $i$-th coordinate of
    $y$ can attain the following values only:
    \begin{align}
      y_i \ &= \ \lambda x_u + (1-\lambda) x_u \ = \ x_u \quad \text{or} \label{eq:yi:uu}\\
      y_i \ &= \ \lambda x_v + (1-\lambda) x_v \ = \ x_v \quad \text{or} \label{eq:yi:vv}\\
      y_i \ &= \ \lambda x_u + (1-\lambda) x_v \ = \ \lambda(x_u-x_v)+x_v \quad
      \text{or} \label{eq:yi:uv} \\
      y_i \ &= \ \lambda x_v + (1-\lambda) x_u \ = \ \lambda(x_v-x_u)+x_u \, . \label{eq:yi:vu}
    \end{align}
    Since $x$ is integral coordinates of types \eqref{eq:yi:uu} and \eqref{eq:yi:vv} are
    integers for arbitrary parameters $\lambda\in(0,1)$.  The coordinates of types
    \eqref{eq:yi:uv} and \eqref{eq:yi:vu} are integral if $\lambda\cdot|x_u-x_v|\in\ZZ$.

    We assumed that $x$ is contained in the $k$-th $\vones$-layer for some
    $k=0,1,\dots,n-1$ and that it is not a core point.  In Example~\ref{exmp:core} it has
    been observed that the core points in these layers are the vertices of a translated
    hypersimplex.  We learned that some coordinate difference $|x_i-x_k|$ must exceed one.
    Since all coefficients are equal to either $x_u$ or $x_v$ it follows that
    $|x_u-x_v|\ge 2$.  We can now set
    \[
    \lambda \ = \ \frac{1}{|x_u-x_v|}
    \]
    in the formula~\eqref{eq:y}.
  \end{enumerate}

  In both cases we obtain an integral point $y$ in the interior of the interval $[x,\gamma
  x]$.  By Lemma~\ref{lem:closer}, such a point is always closer to the center than $x$.
  This shows that there exists a feasible core point in the same layer as $x$.
  Applying Proposition~\ref{prop:mu-trs} with $\mu=\lfloor n/2\rfloor+1$ yields that then
  \emph{each} core point must be feasible.
\end{proof}

Now Algorithm~\ref{algo:reduction_to_layers} can be modified in
Step~\ref{step:reduction_to_layers:feasibility} to check a single core point per layer for
feasibility, provided that the group of automorphisms of the ILP acts at least $(\lfloor
n/2\rfloor+1)$-transitively.  This is our \emph{Core Point
  Algorithm}~\ref{algo:core_point}.

\begin{algorithm}[tb]
  \dontprintsemicolon
  \linesnumbered

  \KwIn{$(A,b)$ such that $\Sym(\ILP(A,b,\vones))$ acts $(\lfloor n/2\rfloor+1)$-transitively on
    standard basis}
  \KwOut{optimal solution of $\ILP(A,b,\vones)$ or ``infeasible''}

  let $z=\zeta\vones$ be a symmetric optimal solution of the LP relaxation $\LP(A,b,\vones)$ \;
  $d\leftarrow \lfloor n\zeta\rfloor-n\lfloor\zeta\rfloor$ \;
  \Repeat{feasible $x$ found or $d<0$}{
    $x\leftarrow(\underbrace{\lfloor \zeta\rfloor+1,\dots,\lfloor
      \zeta\rfloor+1}_{d},\underbrace{\lfloor \zeta\rfloor,\dots,\lfloor
      \zeta\rfloor}_{n-d})$ \;
    \If{$x$ infeasible}{
      $d\leftarrow d-1$ \;
    }
  }
  \Return{$x$ or ``infeasible''}
  \caption{Core point algorithm}
  \label{algo:core_point}
\end{algorithm}

\begin{corollary}\label{cor:core_point}
  For given $A$ and $b$ such that $\Sym(\ILP(A,b,\vones))$ acts $(\lfloor
  n/2\rfloor+1)$-transitively on the standard basis the Core Point
  Algorithm~\ref{algo:core_point} solves the integer linear program $\ILP(A,b,\vones)$ in
  $\order{mn^2}$ time.
\end{corollary}

\begin{proof}
  The correctness follows from Theorems~\ref{thm:reduction_to_layers} and~\ref{thm:core}.
  The main loop of the algorithm is executed at most $n$ times.  In each step the costs
  are dominated by checking one point in $\RR^n$ for feasibility against $m$ linear
  inequalities.
\end{proof}

\begin{remark}
  The linear search in Algorithms~\ref{algo:reduction_to_layers} and~\ref{algo:core_point}
  cannot be substituted by a direct bisectional approach.  The reason is that the set of
  all $k$ in $\{0,1,\dots,\lfloor n\zeta\rfloor-n\lfloor\zeta\rfloor\}$ such that the
  $k$-th $\vones$-layer contains a feasible point is not necessarily (the set of integer
  points of) an interval.
\end{remark}

\section{Finding All Symmetries}
\label{sec:finding}
\noindent 
For the algorithms presented it is never necessary to know the entire group of
automorphisms of $\LP(A,b,c)$ or $\ILP(A,b,c)$.  Generally, any subgroup will do, the
larger the better.  Yet here we would like to discuss the question of how to find
automorphisms of integer linear programs.  From the input data we will construct a labeled
graph $\ILPgraph(A,b,c)$ whose group of labeled automorphisms coincides with
$\Sym(\ILP(A,b,c))$.

Expressing symmetry in optimization via graph automorphisms is not a new idea: the linear
automorphism group of a polytope and of a linear program can be obtained by computing the
automorphism group of a certain graph as described by Bremner, Dutour Sikiri{\'c}, and
Sch\"urmann~\cite{PolyRepConv}. The combinatorial automorphisms of a polytope are the
(labeled) graph automorphisms of the bipartite graph encoded by the
vertex-edge-incidences.  This directly follows from the fact that the face lattice of a
polytope is atomic and coatomic; see Kaibel and Schwartz~\cite{KaibelSchwartz03}.  Liberti
studies automorphisms of optimization problems which are more general than integer linear
programs~\cite{Liberti10}.  His approach, however, deals with expression trees obtained
from a specific encoding of the optimization problem.  None of these concepts seems to be
directly related to the kind of symmetry studied here.  An idea similar to ours, however,
has been applied by Berthold and Pfetsch~\cite{BertholdPfetsch09} to find symmetries of
$0/1$-ILPs.

The complexity status of the graph isomorphism problem is notoriously open.  While the
known algorithms for determining the automorphism group of a graph require exponential
time, there exist software packages, for instance, \texttt{nauty}~\cite{nauty} or
\texttt{SymPol}~\cite{SymPol}, that can solve this problem very well in practice.

For a given matrix $A\in\RR^{m\times n}$, right hand side $b\in\RR^m$, and objective
function $c\in\RR^n$ we will now associate two undirected simple graphs, the \emph{ILP
  graph} $\ILPgraph(A,b,c)$, and the \emph{restricted ILP graph} $\ILPgraph'(A,b,c)$.  For
the sake of a simplified exposition we start out by describing the restricted ILP graph.
Throughout we assume that the rows of the extended matrix $(A|b)$ are normalized as
described in Remark~\ref{rem:scale}.  We have one node $\alpha_{ij}$ for each position in
the matrix $A$, one node $\rho_i$ for each row, and one node $\zeta_j$ for each column,
that is, $(i,j)\in [m]\times[n]$, where $[n]=\{1,2,\dots,n\}$.  Further, we have one node
$\kappa_u$ for each distinct coefficient $u$ in the matrix~$A$, one node $\lambda_v$ for
each distinct coefficient $v$ of $b$, and one node $\mu_w$ for each distinct coefficient
$w$ of $c$.  This gives a total of $mn+m+n+n_A+n_b+n_c$ nodes, where $n_A$, $n_b$, and
$n_c$ denotes the respective number of different entries in $A$, $b$, and $c$.  The nodes
receive labels in the following way: all positions share the same label, the rows receive
a second, and the columns a third label.  Each node corresponding to one of the
coefficients receives an individual label.  This way we arrive at $n_A+n_b+n_c+3$ labels
altogether.  The edges of $\ILPgraph'(A,b,c)$ are defined as follows: the node
$\alpha_{ij}$ is adjacent to $\rho_i$ and $\zeta_j$ as well as to the coefficient node
which represents the coefficient $a_{ij}$ of the matrix $A$.  Moreover, the row node
$\rho_i$ is adjacent to the node $\lambda_{b_i}$, and the node $\zeta_j$ is adjacent to
the node $\mu_{c_j}$.  This totals to $3mn+m+n$ edges.

\begin{example}\label{ex:symILP}
  The reduced ILP graph of the integer linear program
  \begin{equation}\label{eq:symILP}
    \begin{array}{lrcrcrcl}
      \max & x_1&+&x_2&+&x_3&&\\
      \textrm{s.t.}&x_1&+&2x_2&&&\leq&3\\
      &&&x_2&+&2x_3&\leq&3\\
      &2x_1&&&+&x_3&\leq&3 \, , \quad x_i\in\ZZ
    \end{array}
  \end{equation}
  is shown in Figure~\ref{fig:symILP}.

  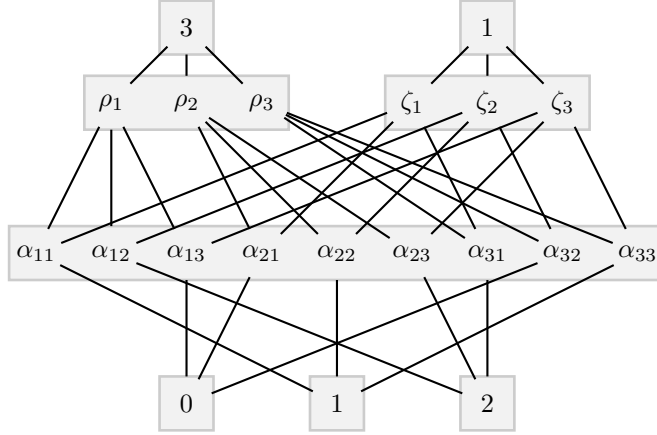
\begin{figure}[htb]
\begin{tikzpicture}
\pgfdeclarelayer{background}
\pgfdeclarelayer{foreground}
\pgfsetlayers{background,main,foreground}
\tikzstyle{vertex}=[circle,minimum size=20pt,inner sep=0pt]
\tikzstyle{partition}=[draw=black!20,rectangle,fill=black!5,line width=1pt,inner sep=0pt]
\tikzstyle{edge} = [draw,thick,-]

    \foreach \pos/\name/\nodeName in {{(2,6)/3/lambda3}, {(6,6)/1/mu1},
      {(2,1)/0/kappa0}, {(4,1)/1/kappa1}, {(6,1)/2/kappa2},
      {(1,5)/\rho_1/rho1}, {(2,5)/\rho_2/rho2}, {(3,5)/\rho_3/rho3},
      {(5,5)/\zeta_1/zeta1}, {(6,5)/\zeta_2/zeta2}, {(7,5)/\zeta_3/zeta3},
      {(0,3)/\alpha_{11}/alpha11}, {(1,3)/\alpha_{12}/alpha12},
      {(2,3)/\alpha_{13}/alpha13}, {(3,3)/\alpha_{21}/alpha21},
      {(4,3)/\alpha_{22}/alpha22}, {(5,3)/\alpha_{23}/alpha23},
      {(6,3)/\alpha_{31}/alpha31}, {(7,3)/\alpha_{32}/alpha32},
      {(8,3)/\alpha_{33}/alpha33}} \node[vertex] (\nodeName) at \pos
             {$\name$};
\begin{pgfonlayer}{background}
  \node[partition,fit=(lambda3)]{}; \node[partition,fit=(mu1)]{};
  \node[partition,fit=(kappa0)]{}; \node[partition,fit=(kappa1)]{};
  \node[partition,fit=(kappa2)]{};
  \node[partition,fit=(rho1) (rho2) (rho3)]{}; 
  \node[partition,fit=(zeta1) (zeta2) (zeta3)]{};
  \node[partition,fit=(alpha11) (alpha12) (alpha13) (alpha21) (alpha22)
    (alpha23) (alpha31) (alpha32) (alpha33)]{};
\end{pgfonlayer}
    \foreach \source/ \dest in {rho1/alpha11, rho1/alpha12,
      rho1/alpha13, rho2/alpha21, rho2/alpha22, rho2/alpha23,
      rho3/alpha31, rho3/alpha32, rho3/alpha33, zeta1/alpha11,
      zeta2/alpha12, zeta3/alpha13, zeta1/alpha21, zeta2/alpha22,
      zeta3/alpha23, zeta1/alpha31, zeta2/alpha32, zeta3/alpha33,
      kappa1/alpha11, kappa2/alpha12, kappa0/alpha13, kappa0/alpha21,
      kappa1/alpha22, kappa2/alpha23, kappa2/alpha31, kappa0/alpha32,
      kappa1/alpha33, zeta1/mu1, zeta2/mu1, zeta3/mu1, rho1/lambda3,
      rho2/lambda3, rho3/lambda3} \path[edge] (\source) -- node {}
    (\dest);
\end{tikzpicture}
    \caption{The reduced ILP graph for \eqref{eq:symILP}.}
    \label{fig:symILP}
  \end{figure}
\end{example}

Let $\gamma$ be an automorphism of $\ILPgraph'(A,b,c)$ which respects all node labels.  Since
the common label of the column nodes is preserved $\gamma$ induces a column permutation
$\psi_\gamma\in\Sym(n)$.  Now $\psi_\gamma$ acts on the standard basis
$\{e_1,e_2,\dots,e_n\}$, and by linear extension we obtain a linear transformation which
we denote $\psi_\gamma^*$.

\begin{lemma}\label{lem:phi_psi}
  The linear transformation $\psi_\gamma^*$ is a symmetry of $\ILP(A,b,c)$.
\end{lemma}

\begin{proof}
  As above let $\gamma$ be a labeled automorphism of $\ILPgraph'=\ILPgraph'(A,b,c)$ with
  induced column permutation $\psi=\psi_\gamma$ and linear transformation
  $\psi^*\in\Sym(n)\le\GL_n\RR$.  As for the column nodes the graph automorphism $\gamma$
  also induces a permutation $\phi\in\Sym(m)$ of the row nodes of $\ILPgraph'$.  The
  position nodes $\alpha_{ij}$ form a label class of their own, and so they are permuted
  by $\gamma$ as well.  Since each position node is adjacent to precisely one row and one
  column node we infer that $\gamma(\alpha_{ij})=\alpha_{\phi(i),\psi(j)}$.  Each position
  node is adjacent to precisely one matrix coefficient node, each of which forms a
  singleton label class.  This implies that the coefficient $a_{ij}$ corresponding to the
  node $\alpha_{ij}$ is the same as the coefficient $a_{\phi(i),\psi(j)}$.  Likewise we
  obtain $b_i=b_{\phi(i)}$ and $c_j=c_{\psi(j)}$.  This means that $\psi_\gamma^*$ is a
  symmetry of $\ILP(A,b,c)$.
\end{proof}

\begin{proposition}\label{prop:reduced_iso}
  The map $\gamma\mapsto\psi_\gamma^*$ is an isomorphism from the group of labeled
  automorphisms of the graph $\ILPgraph'(A,b,c)$ to the group
  $\Sym(\ILP(A,b,c))\cap\Sym(n)$.
\end{proposition}

\begin{proof}
  We describe the inverse map. To this end let
  $\sigma$ be a symmetry of $\ILP(A,b,c)$ which acts on the standard basis of $\RR^n$.
  Hence $\sigma$ induces a permutation $\phi$ of the rows of the extended matrix $(A|b)$
  and a permutation $\psi$ of the columns of $A$.  It is obvious how $\phi$ and $\psi$
  induce permutations of the row nodes and of the column nodes of $\ILPgraph'$.  By the
  same reasoning as in the proof of Lemma~\ref{lem:phi_psi} the pair $(\phi,\psi)$
  uniquely extends to a labeled graph automorphism $\gamma(\sigma)$ of the reduced ILP
  graph.

  We omit the straightforward proofs that the equations $\gamma(\psi_\gamma^*)=\gamma$ and
  $\psi_{\gamma(\sigma)}^*=\sigma$ both hold.  From these it follows that the map
  $\gamma\mapsto\psi_\gamma^*$ is bijective.  In both groups the multiplications are given
  by concatenations of maps.  A direct computation yields
  $\psi_{\gamma_1\gamma_2}^*=\psi_{\gamma_1}^*\psi_{\gamma_2}^*$; all maps are acting on
  the left.  Hence the group structures are preserved.
\end{proof}

We now explain how the full ILP graph $\ILPgraph(A,b,c)$ differs from the restricted ILP
graph $\ILPgraph'(A,b,c)$.  The key to the construction of $\ILPgraph'(A,b,c)$ is the map
$\gamma\mapsto\psi_\gamma^*$ yielding a linear transformation which acts as a permutation
of the standard basis of $\RR^n$.  In order to allow for signed permutations certain nodes
have to be duplicated: each column node $\zeta_j$ in $\ILPgraph'(A,b,c)$ gets a \emph{twin
  node} $\hat\zeta_j$ in $\ILPgraph(A,b,c)$, each matrix coefficient node $\alpha_{ij}$
corresponding to a non-zero coefficient gets a twin node $\hat\alpha_{ij}$.  Moreover, we
add further nodes representing negatives of non-zero coefficients in the matrix $A$ and
the objective function $c$ unless nodes with these labels already exist.  This way
$\ILP(A,b,c)$ has less than twice as many nodes as $\ILP'(A,b,c)$; it is always strictly
less as the nodes corresponding to the coefficients in $b$ are never duplicated.  We also
add edges such that first $\hat\alpha_{ij}$ is adjacent to $\rho_i$ and $\hat\zeta_j$ for
all $i$ and $j$, second $\hat\zeta_j$ is adjacent to $\mu_{-c_j}$, third $\hat\alpha_{ij}$
is adjacent to $\kappa_{-a_{ij}}$, and, finally, the twins are matched up: $\alpha_{ij}$
is adjacent to $\hat\alpha_{ij}$ and $\hat\zeta_j$ is adjacent to $\zeta_j$.  The labeling
is extended in a way such that twins share the same label; the nodes newly introduced for
negatives of coefficients receive new singleton labels.

Each labeled graph automorphism of $\ILPgraph'(A,b,c)$ uniquely extends to a labeled graph
automorphism of $\ILPgraph(A,b,c)$, but the automorphism group of the non-reduced ILP
graph is larger, in general.  We have the following result.

\begin{theorem}
  The group of labeled graph automorphisms of $\ILPgraph(A,b,c)$ is isomorphic to the
  group of symmetries of $\ILP(A,b,c)$.
\end{theorem}

\begin{proof}
  One can follow the strategy in the proof of Proposition~\ref{prop:reduced_iso}.  We know
  that a labeled graph automorphism of $\ILPgraph'(A,b,c)$ encodes a symmetry of
  $\ILP(A,b,c)$ which permutes the set $\{e_1,e_2,\dots,e_n\}$.  Now a labeled graph
  automorphism of $\ILPgraph(A,b,c)$ may map a column node $\zeta_j$ to some node
  $\hat\zeta_k$.  But then it follows that $\hat\zeta_j$ is mapped to $\zeta_k$ since
  $\hat\zeta_j$ is the only column node adjacent to $\zeta_j$, and $\zeta_k$ is the only
  column node adjacent to $\hat\zeta_k$.  This shows that the permutation of the column
  nodes can be extended to a linear transformation.  As in the proof of
  Proposition~\ref{prop:reduced_iso} one can show that this linear transformation is a
  symmetry of the integer linear program.  Conversely, each such symmetry acts like a
  signed permutation on the signed standard basis and yields a labeled isomorphism of the
  graph $\ILPgraph(A,b,c)$.
\end{proof}

Roughly speaking, a class $\mathcal{C}$ of graphs is \emph{graph isomorphism complete} if
the problem of deciding isomorphy for any two graphs in $\mathcal{C}$ is as difficult as
for general graphs, up to a polynomial time transformation.  For a precise definition, for
instance, see the monograph~\cite{KoeblerSchoeningToran93}.  The next result is not only
of theoretical interest.  To the contrary, for practical applications it can be read as:
finding the symmetries of an integer linear program via reducing to automorphisms of
suitable (labeled) graphs, is the right thing to do.

\begin{theorem}
  The classes of ILP graphs and reduced ILP graphs are both graph isomorphism complete.
\end{theorem}

\begin{proof}
  We only prove that the class of reduced ILP graphs is graph isomorphism complete.  It is
  known that the class of bipartite graphs is graph isomorphism complete.  Hence it
  suffices to encode an arbitrary bipartite graph as a reduced ILP graph, which is not too
  large.

  Let $G=(V,E)$ be an undirected bipartite graph with $m+n$ nodes $V=U\cup
  W=\{u_1,\dots,u_m\}\cup\{w_1,\dots,w_n\}$.  As our matrix $A_G=(a_{ij})\in\RR^{m\times
    n}$ we take the bipartite adjacency matrix of $G$, that is,
  \[
  a_{ij} \ = \
  \begin{cases}
    1 & \text{if $\{u_i,w_j\}\in E$}\\
    0 & \text{otherwise \, .}
  \end{cases}
  \]
  For a second bipartite graph $G'$ it is easy to see that the reduced ILP graph of
  $\ILP(A_G,\vones,\vones)$ is isomorphic to the reduced ILP graph of
  $\ILP(A_{G'},\vones,\vones)$ if and only if $G$ is isomorphic to $G'$.
\end{proof}

\begin{remark}\label{rem:GLnZ}
  Rehn investigates arbitrary automorphisms of the integer lattice $\ZZ^n$ in the context
  of polyhedral geometry~\cite{Rehn10}.  In particular, his modification of a backtracking
  algorithm of Plesken and Souvignier~\cite{PleskenSouvignier97} allows to obtain matrix
  generators of the group of symmetries.  For practical applications this should be
  superior to our approach via graph automorphisms if the number $m$ of constraints is
  much larger than the dimension $n$.
\end{remark}

\section{Hypertruncated Cubes}
\label{sec:hypertruncated}
\noindent
In this section we will construct a specific class of highly symmetric convex polytopes
among which one can find examples of rather high Gomory-Chv\'atal rank.  The motivation
for this construction is rooted in the systematic study of Gomory-Chv\'atal cuts and
cutting-plane proof systems.  Pokutta and Stauffer~\cite{PokuttaStauffer10} propose a new
method for computing lower bounds on the Gomory-Chv\'atal rank for polytopes contained in
the $0/1$-cube, and the polytopes constructed here provide examples which asymptotically
almost attain the bounds obtained.  The subsequent section on computational experiments
also contains results about these polytopes.

Our construction starts out with the unit cube $C=[0,1]^n$.  Intersecting $C$ with the
hyperplane defined by $\sum x_i=r$ for $r\in\{2,3,\dots,n-1\}$ gives the hypersimplex
$\Delta(r,n)$ which already appeared in Example~\ref{exmp:core}.  Here we are interested in
the \emph{$(r,n)$-truncated cube} $C' = \smallSetOf{x\in [0,1]^n}{\sum x_i\le r}$.  We
make one more modification to the cube $C$ by defining the polytope
\begin{equation}\label{eq:hypertruncated}
  C'' \ = \ \conv(C'\cup\lambda\vones) \quad \text{for $\lambda>r/n$}
\end{equation}
which we call the \emph{$(r,n;\lambda)$-hypertruncated cube}.  Notice that the full group
$\Sym(n)$ acts on the cube $C$ as well as on the truncated cube $C'$ as well as on the
hypertruncated cube $C''$.  Hence our algorithms above can be applied.  Our next goal is
to describe the vertices and the facets of $C''$.

\begin{proposition}
  Let $n\ge 2$, $r\in\{2,3,\dots,n-1\}$, and $r/n<\lambda<1$.  The vertices of the
  $(r,n;\lambda)$-hypertruncated cube $C''$ are
  \[
  e_S \quad \text{for all $S\subset[n]$ with $\#S\le r$} \qquad \text{and} \qquad
  \lambda\vones \, .
  \]
\end{proposition}

\begin{proof}
  The points $e_S$, for $S\subset[n]$ and $\#S\le r$, are the vertices of the
  $(r,n)$-truncated cube $C'$.  They are also vertices of $C''$.  Since $n\lambda$ exceeds
  $r$, the hyperplane $\sum x_i=n\lambda$ does not separate $C''$, and its intersection
  with $C''$ only contains the point $\lambda\vones$.  Hence the latter is a vertex, too.
  Looking at the defining Equation~\eqref{eq:hypertruncated} shows that there cannot be
  any other vertices.
\end{proof}

Of course, the vertices determine the facets completely.  In this case, it is particularly
easy to read off the facets of $C''$ by looking at the facets of $C'$ and analyzing what
changes in case the point $\lambda\vones$ is added as a generator.  This proves the claim
in~\cite[Remark~3.3]{PokuttaStauffer10}.

\begin{proposition}\label{prop:hypertruncated_facets}
  Let $n\ge 2$, $r\in\{2,3,\dots,n-1\}$, and $\lambda>r/n$.  The facets of the
  $(r,n;\lambda)$-hypertruncated cube $C''$ are
  \begin{align}
    & x_i \ \ge \ 0 \ , \quad x_i \ \le \ 1 \label{eq:hypertruncated_facets:cube}\\
    & \bigl(1-n+\frac{r}{\lambda}\bigr) x_i + \sum_{k\ne i} x_k \ \ge \ r  \label{eq:hypertruncated_facets:deletion}\\
    & \bigl(1-r+\lambda(n-1)\bigr) x_i + (1-\lambda) \sum_{k\ne i} x_k \ \le \
    \lambda(n-r) \label{eq:hypertruncated_facets:contraction}
  \end{align}
  for $i\in[n]$.  In particular, $C''$ has precisely $4n$ facets.
\end{proposition}

\begin{proof}
  The facets of type~\eqref{eq:hypertruncated_facets:cube} are the facets of the unit cube
  $C$.  Together with the truncating inequality $\sum x_i\le r$ they also form the facets
  of the truncated cube~$C'$.  The remaining facets of $C''$ are the facets through the
  vertex~$\lambda\vones$.  Each of them is the convex hull of $\lambda\vones$ and a ridge
  of $C'$ contained in the truncating facet.  A \emph{ridge} is a face of codimension~$2$,
  that is, a facet of a facet.  As pointed out above the truncating facet is the
  hypersimplex $\Delta(r,n)$.  Its facets arise from the intersection with the cube
  facets.  A hypersimplex facet of type $\smallSetOf{x\in\Delta(r,n)}{x_i=0}$ is a
  \emph{deletion facet}, and a hypersimplex facet of type
  $\smallSetOf{x\in\Delta(r,n)}{x_i=1}$ is a \emph{contraction facet}.  The $n-1$ points
  $r e_k$ for $k\ne i$ span an $(n-2)$-dimensional affine subspace $\cA$ containing the
  $i$-th deletion facet.  However, these points are not contained in $\Delta(r,n)$.
  Looking for an affine hyperplane containing $\cA$ and $\lambda\vones$ results in a
  rank-$1$ system of linear equations.  This way we obtain the $n$ linear inequalities of
  type~\eqref{eq:hypertruncated_facets:deletion}.  Similarly, the affine span of a
  contraction ridge is generated by the $n-1$ points $e_i+(r-1)e_k$ for $k\ne i$.  Via the
  same approach we arrive at the $n$ linear inequalities of
  type~\eqref{eq:hypertruncated_facets:contraction}.
\end{proof}

\begin{remark}\label{rem:hypertruncated:parameters}
  Pokutta and Stauffer~\cite{PokuttaStauffer10} show that the Gomory-Chv\'atal ranks of
  the $(r,n;\lambda)$-hyper\-trun\-ca\-ted cubes for $r=\lfloor n/e\rfloor$, where
  $e=2.7172\ldots$ is Euler's constant, and $\lambda=(m-1)/m$ approach $n/e - o(1)$ as
  $m\in\NN$ goes to infinity.  In our experiments below we look at the case $r=\lfloor
  n/e\rfloor$ and $\lambda=1/2$, that is, $m=2$.
\end{remark}

\section{Computational Results}
\label{sec:experiments}
\noindent
The following experiments were carried out on an Intel(R) Core(TM) i7\ 920\ 2.67GHz
machine, with 12GB of main memory, running Ubuntu 10.04 (Lucid Lynx).  The performance of
each core is estimated at 5346.16 bogomips each.  All our tests were run single-threaded.

The goal of the experiments is to compare the performances of a conventional
branch-and-cut approach (with automated symmetry detection) and the Core-Point-Algorithm
\ref{algo:core_point} on highly symmetric integer linear programs. For the test of the
conventional branch-and-cut method we used \texttt{CPLEX}, Version~12.1.0, while the
Core-Point-Algorithm was implemented and tested in \texttt{polymake},
Version~2.9.9~\cite{polymake}.  As a major difference \texttt{polymake} employs exact
rational arithmetic (via \texttt{GMP}), while \texttt{CPLEX} uses floating-point
arithmetic.  It should be stressed that \texttt{CPLEX} can detect if the symmetry group of
an integer linear program contains the full symmetric group acting on the standard basis
of $\RR^n$, and this is exploited in its algorithms.  For input in this category (which
includes all our examples below), it is thus quite a challenge to beat \texttt{CPLEX}.

\subsection{Hypertruncated Cubes}

We tested our algorithms on the $(\lfloor n/e\rfloor,n;1/2)$-hyper\-trunca\-ted cubes; see
Remark~\ref{rem:hypertruncated:parameters}.  In this case we have only $4n$ linear
inequalities from Proposition~\ref{prop:hypertruncated_facets} as input.  Each coefficient
is small, and computationally accuracy (for floating-point computations) or coefficient
growth (for exact arithmetic) is not an issue here.  This benign input can be dealt with
easily up to high dimensions.  Table~\ref{tab:HypertruncatedCube} lists the timings for
\texttt{CPLEX}' Branch-and-Cut and \texttt{polymake}'s Core Point Algorithm.  The timings
required to obtain the solution of the linear relaxation are given separately for both
systems.

The fact that \texttt{polymake} takes more time is due to the overhead induced by the
\texttt{GMP} exact rational arithmetic.  Since coefficient growth does not occur the
overhead versus floating-point arithmetic can be estimated to be constant.  Hence the
roughly quadratic overhead (in dependence of $n$) versus the \texttt{CPLEX} result is a
consequence of the total algorithmic complexity of $\order{mn^2}$ from
Corollary~\ref{cor:core_point}.  Altogether both solvers behave pretty well for these
kinds of examples.

\begin{table}[ht]
\caption{Hypertruncated Cubes}
  \label{tab:HypertruncatedCube}
  \renewcommand{\arraystretch}{0.9}
  \begin{tabular*}{\linewidth}{@{\extracolsep{\fill}}rrrrr@{}}
    \toprule
    &  \multicolumn{2}{c}{\texttt{CPLEX}} &  \multicolumn{2}{c}{\texttt{polymake}} \\
    \multicolumn{1}{@{}r}{$d$} & \multicolumn{1}{r}{time LP (s)} & \multicolumn{1}{r}{time
    IP (s)} & \multicolumn{1}{r}{time LP (s)} & \multicolumn{1}{r@{}}{time IP (s)}\\
    \midrule
    100 & 0.00 & 0.07 & 0.01 & 0.08\\
    200 & 0.07 & 0.29 & 0.02 & 0.57\\
    300 & 0.19 & 0.73 & 0.03 & 1.88\\
    400 & 0.41 & 1.58 & 0.06 & 4.26\\
    500 & 0.90 & 2.99 & 0.10 & 8.39\\
    600 & 1.52 & 4.80 & 0.14 & 14.30\\
    700 & 2.21 & 7.01 & 0.18 & 22.29\\
    800 & 3.44 & 11.59 & 0.24 & 32.96\\
    900 & 5.17 & 16.37 & 0.31 & 47.11\\
    1000 & 6.77 & 21.66 & 0.38 & 65.11\\
    1100 & 8.25 & 26.55 & 0.47 & 85.52\\
    1200 & 11.75 & 35.47 & 0.55 & 111.59\\
    1300 & 14.05 & 45.63 & 0.65 & 142.76\\
    1400 & 18.96 & 57.23 & 0.74 & 175.23\\
    1500 & 23.42 & 73.50 & 0.85 & 217.19\\
    1600 & 28.11 & 78.13 & 0.98 & 263.21\\
    1700 & 32.07 & 97.23 & 1.11 & 315.38\\
    1800 & 41.82 & 128.88 & 1.25 & 374.63\\
    1900 & 44.68 & 137.22 & 1.40 & 444.48\\
    2000 & 50.39 & 154.35 & 1.54 & 511.59\\
    \bottomrule
  \end{tabular*}
\end{table}

An industry strength solver as \texttt{CPLEX} comes with a number of bolts and whistles
which allow to tune its behavior in many ways.  For the hypertruncated cubes this does not
play any role.  Since no parallel implementation of the Core Point Algorithm is available
(yet) we set the number of \texttt{CPLEX}' parallel threads to one.

\subsection{Wild Input}

One way to produce symmetric input for (integer) linear optimization algorithms is by
brute force: One can take any system $Ax\le b$ of linear inequalities and let the full
group $\Sym(n)$ act.  This way each original inequality may give up to $n!$ inequalities
in the resulting \emph{symmetrized} system.  The symmetrized system of inequalities is
$\Sym(n)$-invariant by construction.  In order to produce input to our algorithms which is
less well behaved than the hypertruncated cubes studied above we will apply this procedure
to a special class of polytopes, which can be considered ``wild''.  We aim at symmetric
polytopes with many facets whose coordinates are not so nice, but still somewhat under
control.

The first building block of our construction is the regular hexagon $H$ whose vertices are
at distance $56/6$ from the origin, that is,
\[
H \ = \ \conv \SetOf{\frac{56}{6} e^{k\pi i/6}}{k=0,1,\dots,5} \, .
\]
Notice that only in the formula above the letter `$i$' denotes the imaginary unit, and we
identify the complex numbers with $\RR^2$.  The coordinates of $H$ are irrational;
however, the subsequent steps in the construction are chosen such that we will arrive at a
rational polytope in the end.  The second item is the regular cross polytope scaled by
$73/10$, that is,
\[
C(d) \ = \ \conv\SetOf{\pm \frac{73}{10} e_i}{i\in[d]} \, .
\]
Finally, we consider the \emph{join} $P*Q$ of two polytopes $P\subset\RR^\delta$ and
$Q\subset\RR^\epsilon$, which is defined as
\[
\begin{split}
  P*Q \ = \ \conv\bigl(&\SetOf{(x,0,1)\in\RR^\delta\times\RR^\epsilon\times\RR}{x\in P}\\
  &\cup \SetOf{(0,y,-1)\in\RR^\delta\times\RR^\epsilon\times\RR}{y\in Q}\bigr) \, .
\end{split}
\]
If $P$ and $Q$ are full-dimensional polytopes with $\mu$ and $\nu$ vertices,
respectively, the join $P*Q$ has dimension $\delta+\epsilon+1$ and $\mu+\nu$ vertices.
For the combinatorics of $P*Q$ the exact values for the $(\delta+\epsilon+1)$st coordinate
are inessential, as long as they are distinct.  We now replace the ``$-1$'' for the
second factor by $-11/12$ to obtain the \emph{distorted join}
\[
\begin{split}
J(d) \ = \ \conv\bigl(&\SetOf{(x,0,1)\in\RR^2\times\RR^d\times\RR}{x\in H}\\
  &\cup \SetOf{(0,y,-\frac{11}{12})\in\RR^2\times\RR^d\times\RR}{y\in C(d)}\bigr) \,
\end{split}
\]
of the hexagon $H$ with the cross polytope $C(d)$.  This polytope is further modified in
two steps: First we perturb by rounding the (rational and irrational) coordinates to three
decimal places and treating these as exact rational numbers.  Since the polytopes $H$,
$C(d)$, and $J(d)$ are simplicial this perturbation does not change the combinatorial
types.  Secondly, we symmetrize the polytope by letting the group $\Sym(d+3)$ act on the
facets of the perturbed polytope.  The resulting inequalities form the input of our second
class of experiments.

The parameters $56/6$, $73/10$, and $11/12$ which occur in the construction are chosen,
more or less, at random.  They do not have a specific meaning.  We refrain from further
investigating these symmetrized distorted joins and the geometry of lattice points inside.
This would be tedious and at the same time irrelevant for our purposes.

\begin{table}
\caption{Symmetrized distorted joins of a hexagon with cross-polytopes}
  \label{tab:joinNgonCross}
  \renewcommand{\arraystretch}{0.9}
  \begin{tabular*}{\linewidth}{@{\extracolsep{\fill}}rrrrr@{}}
    \toprule
    &  \multicolumn{2}{c}{\texttt{CPLEX}} &  \multicolumn{2}{c}{\texttt{polymake}} \\
    \multicolumn{1}{@{}r}{$d$} & \multicolumn{1}{r}{time LP (s)} & \multicolumn{1}{r}{time
    IP (s)} & \multicolumn{1}{r}{time
    LP (s)} & \multicolumn{1}{r@{}}{time IP (s)}\\
    \midrule
    3 & 0.00 & 0.01 & 0.00 & 0.00\\
    4 & 0.00 & 0.06 & 0.01 & 0.00\\
    5 & 0.00 & 0.17 & 0.01 & 0.02\\
    6 & 0.05 & 0.74 & 0.04 & 0.04\\
    7 & 0.13 & 2.71 & 0.09 & 0.13\\
    8 & 0.62 & 10.15 & 0.24 & 0.38\\
    9 & 2.08 & 42.06 & 0.69 & 1.03\\
    10 & 8.02 & 135.51 & 1.86 & 2.89\\
    \bottomrule
  \end{tabular*}
\end{table}

The interesting fact is that we get symmetric polytopes which are somewhat complicated,
because they have lots of inequalities: for instance, yielding $885,\!768$ inequalities
for $d=10$.  As a consequence \texttt{CPLEX} cannot deal with these examples in a fully
automated way.  The best parameter settings that we found were
\begin{center}
  \begin{tabular*}{.7\linewidth}{lr}
    parallel thread count: & 1\\
    presolve indicator: & no\\
    feas.\ pump heuristic: & -1\\
    RINS heuristic: & -1\\
    MIP optimization emph.: &2
  \end{tabular*}
\end{center}
But even with these adjustments our implementation outperforms \texttt{CPLEX} by a large
margin; see Table~\ref{tab:joinNgonCross}.  This holds in spite of the fact that
\texttt{polymake} computes with exact rational numbers throughout.

\bibliographystyle{amsplain}
\bibliography{main}

\end{document}